\documentclass[11pt]{article}

\usepackage{amssymb,latexsym}
\usepackage{epsfig}
\usepackage{eufrak}
\usepackage{amsmath}
\usepackage{mathrsfs}
\usepackage{color}

\newtheorem{theorem}{Theorem}
\newtheorem{lemma}{Lemma}

\newcommand{\be}{\begin{equation}}
\newcommand{\ee}{\end{equation}}
\newcommand{\bee}{\begin{eqnarray*}}
\newcommand{\eee}{\end{eqnarray*}}
\newcommand{\bel}{\begin{eqnarray}}
\newcommand{\eel}{\end{eqnarray}}
\newcommand{\bec}{\begin{cases}}
\newcommand{\eec}{\end{cases}}
\newcommand{\bem}{\begin{bmatrix}}
\newcommand{\eem}{\end{bmatrix}}

\newcommand{\la}{\label}
\newcommand{\li}{\left}
\newcommand{\ri}{\right}

\newcommand{\vep}{\varepsilon}
\newcommand{\lm}{\lambda}

\newcommand{\de}{\delta}

\newcommand{\se}{\theta}

\newcommand{\f}{\frac}

\newcommand{\cd}{\cdots}

\newcommand{\bb}{\mathbb}

\newcommand{\wh}{\widehat}

\newcommand{\mrm}{\mathrm}
\newcommand{\bs}{\boldsymbol}

\newcommand{\iy}{\infty}

\newcommand{\pa}{\partial}

\newcommand{\bed}{\begin{description}}
\newcommand{\eed}{\end{description}}
\newcommand{\bei}{\begin{itemize}}
\newcommand{\eei}{\end{itemize}}
\newcommand{\ben}{\begin{enumerate}}
\newcommand{\een}{\end{enumerate}}
\newcommand{\bib}{\bibitem}
\newcommand{\beL}{\begin{lemma}}
\newcommand{\eeL}{\end{lemma}}
\newcommand{\beT}{\begin{theorem}}
\newcommand{\eeT}{\end{theorem}}
\newcommand{\sect}{\section}

\newcommand{\bpf}{\begin{pf}}
\newcommand{\epf}{\end{pf}}
\newcommand{\bsk}{\bigskip}

\newcommand{\pfbox}{\hfill\mbox{$\Box$}}

\newenvironment{pf}{\paragraph*{Proof{\rm.}}}{\pfbox\bigskip}

\begin{document}

\title{{\bf A Simple Sample Size Formula for Estimating Means of Poisson Random Variables}
\thanks{The author is currently with Department of Electrical Engineering,
Louisiana State University at Baton Rouge, LA 70803, USA, and Department of Electrical Engineering, Southern University and A\&M College, Baton
Rouge, LA 70813, USA; Email: chenxinjia@gmail.com}}

\author{Xinjia Chen}

\date{Submitted in April, 2008}

\maketitle

\begin{abstract}

In this paper, we derive an explicit sample size formula based a mixed criterion of absolute and relative errors for estimating means of Poisson
random variables.

\end{abstract}

\section{Sample Size Formula}

It is a frequent problem to estimate the mean value of a Poisson random variable based on sampling.  Specifically, let $X$ be a Poisson random
variable with mean $\bb{E}[ X ] = \lm > 0$, one wishes to estimate $\lm$ as
\[
\wh{\bs{\lm}} = \f{\sum_{i = 1}^n X_i}{n}
\]
where $X_1, \cd, X_n$ are i.i.d. random samples of $X$.  Since $\wh{\bs{\lm}}$ is of random nature, it is important to control the statistical
error of the estimate. For this purpose, we have

\beT Let $\vep_a > 0, \; \vep_r \in (0, 1)$ and $\de \in (0,1)$.  Then
\[
\Pr \li \{  \li | \wh{\bs{\lm}}  - \lm \ri | < \vep_a \; \mrm{or} \; \li | \wh{\bs{\lm}}  - \lm \ri | < \vep_r \lm \ri \} > 1 - \de
\]
provided that \be \la{ccc}
 n > \f{ \vep_r } { \vep_a } \times \f{ \ln \f{2}{\de} } {  (1 + \vep_r) \ln (1 + \vep_r)  - \vep_r }.
\ee \eeT

\bsk

It should be noted that conventional methods for determining sample sizes are based on normal approximation, see \cite{Desu} and the references
therein.  In contrast, Theorem 1 offers a rigorous method for determining sample sizes.  To reduce conservatism, a numerical approach has been
developed by Chen \cite{Chen} which permits exact computation of the minimum sample size.

\sect{Proof of Theorem 1}

We need some preliminary results.

\beL \la{lem1} Let $K$ be a Poisson random variable with mean $\se > 0$. Then, $\Pr \{ K \geq r \} \leq e^{- \se} \li ( \f{ \se e } { r } \ri
)^r$ for any real number $r > \se$ and $\Pr \{ K \leq r \} \leq e^{- \se} \li ( \f{ \se e } { r } \ri )^r$ for any positive real number $r <
\se$. \eeL

\bpf For any real number $r > \se$, using the Chernoff bound \cite{Chernoff}, we have \bee \Pr \{ K \geq r \} & \leq & \inf_{t > 0} \bb{E} \li [
e^{ t (K - r) } \ri ] = \inf_{t
> 0} \sum_{i
= 0}^\iy e^{ t (i - r)} \f{ \se^i } { i! } e^{- \se}\\
& = & \inf_{t > 0} e^{\se e^t } e^{- \se} e^{ - r \; t} \sum_{i = 0}^\iy \f{ (\se e^t)^i } { i! } e^{- \se e^t} =  \inf_{t > 0} e^{- \se} e^{
\se e^t - r \; t}, \eee where the infimum is achieved at $t = \ln \li ( \f{r}{\se} \ri ) > 0$.  For this value of $t$, we have $e^{- \se} e^{
\se e^t - t r } = e^{- \se} \li ( \f{ \se e } { r } \ri )^r$.  It follows that $\Pr \{ K \geq r \} \leq e^{- \se} \li ( \f{ \se e } { r } \ri
)^r$ for any real number $r > \se$.

Similarly, for any real number $r < \se$, we have $\Pr \{ K \leq r \} \leq e^{- \se} \li ( \f{ \se e } { r } \ri )^r$.

\epf

In the sequel, we shall introduce the following function
\[
g(\vep, \lm) = \vep + (\lm + \vep ) \ln \f{\lm}{\lm + \vep}.
\]

\beL \la{lem2}
Let $\lm > \vep > 0$.  Then, $\Pr \li \{ \wh{\bs{\lm}} \leq \lm - \vep \ri \}  \leq  \exp \li ( n \; g(- \vep, \lm) \ri )$ and
$g(- \vep, \lm)$ is monotonically increasing with respect to $\lm \in (\vep, \iy)$. \eeL

\bpf

Letting $K = \sum_{i=1}^n X_i, \; \se = n \lm$ and $r = n (\lm - \vep)$ and applying Lemma \ref{lem1},  for $\lm > \vep > 0$,  we have \[  \Pr
\li \{ \wh{\bs{\lm}} \leq \lm - \vep \ri \} =  \Pr \{ K \leq r \}  \leq  e^{- \se} \li ( \f{ \se e } { r } \ri )^r  =  \exp \li ( n \; g(- \vep,
\lm) \ri ), \] where $g(- \vep, \lm)$ is monotonically increasing with respect to $\lm \in (\vep, \iy)$ because
\[
\f{\pa g(- \vep, \lm) } { \pa \lm  } = - \ln \li ( 1 - \f{\vep}{\lm} \ri ) - \f{\vep}{\lm} >  0
\]
for $\lm > \vep > 0$.

\epf

\beL \la{lem3} Let $\vep > 0$.  Then, $\Pr \li \{ \wh{\bs{\lm}} \geq \lm + \vep \ri \}  \leq  \exp \li ( n \; g(\vep, \lm) \ri )$ and $g( \vep,
\lm)$ is monotonically increasing with respect to $\lm \in (0, \iy)$. \eeL

\bpf

Letting $K = \sum_{i=1}^n X_i, \; \se = n \lm$ and $r = n (\lm + \vep)$ and applying Lemma \ref{lem1},  for $\lm > 0$,  we have
\[ \Pr \li \{ \wh{\bs{\lm}} \geq \lm + \vep \ri \}  =  \Pr \{ K \geq r \} \leq  e^{- \se} \li ( \f{ \se e } { r } \ri )^r \leq \exp \li ( n \;
g(\vep, \lm) \ri ), \]  where $g(\vep, \lm)$ is monotonically increasing with respect to $\lm \in (0, \iy)$ because
\[
\f{\pa g(\vep, \lm) } { \pa \lm  } = - \ln \li ( 1 + \f{\vep}{\lm} \ri ) + \f{\vep}{\lm} > 0.
\]

\epf

\beL \la{lem4} $g (\vep, \lm) > g (- \vep, \lm)$ for $\lm > \vep > 0$. \eeL

\bpf
 Since $g (\vep, \lm) - g (- \vep, \lm) = 0$ for $\vep = 0$ and
\[
\f{ \pa \; [ g (\vep, \lm) - g (- \vep, \lm) ] } { \pa \vep  } = \ln \f{ \lm^2 } { \lm^2 - \vep^2 } > 0
\]
for $\lm > \vep > 0$, we have
\[
g (\vep, \lm) - g (- \vep, \lm)  > 0
\]
for any $\vep \in (0, \lm)$.  Since such arguments hold for arbitrary $\lm > 0$, we can conclude that
\[
g (\vep, \lm) > g (- \vep, \lm)
\]
for $\lm > \vep > 0$. \epf

\beL \la{lem5} Let $0 < \vep < 1$.  Then, $\Pr \li \{ \wh{\bs{\lm}} \leq \lm (1 - \vep) \ri \}  \leq  \exp \li ( n \; g (-\vep \lm, \lm)  \ri )$
and $g (-\vep \lm, \lm)$ is monotonically decreasing with respect to $\lm > 0$. \eeL

\bpf

Letting $K = \sum_{i=1}^n X_i, \; \se = n \lm$ and $r = n \lm (1 - \vep)$ and making use of Lemma \ref{lem1}, for $0 < \vep < 1$, we have
 \[  \Pr \li \{ \wh{\bs{\lm}} \leq \lm (1 - \vep) \ri
\} =  \Pr \{ K \leq r \}  \leq  e^{- \se} \li ( \f{ \se e } { r } \ri )^r   \leq  \exp \li ( n \; g (-\vep \lm, \lm)  \ri ),  \] where
\[
g (-\vep \lm, \lm)  = \li [ - \vep - (1 - \vep ) \ln (1 - \vep) \ri ] \lm,
\]
which is monotonically decreasing with respect to $\lm > 0$, since $- \vep - (1 - \vep ) \ln (1 - \vep) < 0$ for $0 < \vep < 1$.

\epf

\beL \la{lem6} Let $\vep > 0$.  Then, $\Pr \li \{ \wh{\bs{\lm}} \geq \lm (1 + \vep) \ri \}  \leq  \exp \li ( n \; g (\vep \lm, \lm)  \ri )$ and
$g (\vep \lm, \lm)$ is monotonically decreasing with respect to $\lm > 0$. \eeL

\bpf

Letting $K = \sum_{i=1}^n X_i, \; \se = n \lm$ and $r = n \lm (1 + \vep)$ and  making use of Lemma \ref{lem1}, for $\vep > 0$,  we have
 \bee  \Pr \li \{ \wh{\bs{\lm}} \geq \lm (1 + \vep) \ri
\} & \leq & \exp \li ( n \; g (\vep \lm, \lm)  \ri ) \eee where
\[
g (\vep \lm, \lm) = \li [ \vep - (1 + \vep ) \ln (1 + \vep) \ri ] \lm,
\]
which is monotonically decreasing with respect to $\lm > 0$, since $\vep - (1 + \vep ) \ln (1 + \vep) < 0$ for $\vep > 0$.

\epf

\bsk

We are now in a position to prove the theorem.  It suffices to show \[ \Pr \li \{  \li | \wh{\bs{\lm}}  - \lm \ri | \geq \vep_a \; \& \;  \li |
\wh{\bs{\lm}}  - \lm \ri | \geq \vep_r \lm \ri \} < \de
\]
for $n$ satisfying (\ref{ccc}).  It can shown that (\ref{ccc}) is equivalent to \be \la{nnn}
 \exp( n \; g (\vep_a, \vep_a) ) < \f{\de}{2}.
 \ee

We shall consider four cases as follows.

Case (i): $0 < \lm < \vep_a$;

Case (ii): $\lm = \vep_a$;

Case (iii): $\vep_a < \lm \leq \f{\vep_a}{\vep_r}$;

Case (iv): $\lm > \f{\vep_a}{\vep_r}$.

\bsk

In Case (i), we have $\Pr \{ \wh{\bs{\lm}} \leq \lm -  \vep_a \} = 0$ and \bee \Pr \li \{  \li | \wh{\bs{\lm}}  - \lm \ri | \geq \vep_a \; \& \;
\li | \wh{\bs{\lm}}  - \lm \ri | \geq \vep_r \lm
\ri \} & = & \Pr \li \{  \li | \wh{\bs{\lm}}  - \lm \ri | \geq \vep_a  \ri \}\\
& = & \Pr \{ \wh{\bs{\lm}} \leq \lm -  \vep_a \} + \Pr \{ \wh{\bs{\lm}} \geq \lm  + \vep_a \}\\
& = & \Pr \{ \wh{\bs{\lm}} \geq \lm  + \vep_a \}. \eee By Lemma (\ref{lem3}),
\[
\Pr \{ \wh{\bs{\lm}} \geq \lm  + \vep_a \} \leq \exp( n \; g (\vep_a, \lm) ) \leq \exp( n \; g (\vep_a, \vep_a) ) < \f{\de}{2}.
\]
Hence,
\[
\Pr \li \{  \li | \wh{\bs{\lm}}  - \lm \ri | \geq \vep_a \; \& \; \li | \wh{\bs{\lm}}  - \lm \ri | \geq \vep_r \lm \ri \} < \f{\de}{2} < \de.
\]

\bsk

In Case (ii), we have $\Pr \{ \wh{\bs{\lm}} \leq \lm -  \vep_a \} = \Pr \{ \wh{\bs{\lm}} = 0 \}$ and \bee \Pr \li \{  \li | \wh{\bs{\lm}} - \lm
\ri | \geq \vep_a \; \& \; \li | \wh{\bs{\lm}}  - \lm \ri | \geq \vep_r \lm
\ri \} & = & \Pr \li \{  \li | \wh{\bs{\lm}}  - \lm \ri | \geq \vep_a  \ri \}\\
& = & \Pr \{ \wh{\bs{\lm}} \leq \lm -  \vep_a \} + \Pr \{ \wh{\bs{\lm}} \geq \lm  + \vep_a \}\\
& = & \Pr \{ \wh{\bs{\lm}} = 0 \} + \Pr \{ \wh{\bs{\lm}} \geq \lm  + \vep_a \}. \eee Noting that $\ln 2 < 1$, we can show that $- \vep_a < g
(\vep_a, \vep_a)$ and hence  \bee
 \Pr \{ \wh{\bs{\lm}} = 0 \} & = & \Pr \{ X_i = 0, \; i = 1, \cd, n \}\\
 & = & [ \Pr \{ X = 0 \} ]^n \\
 & = & e^{- n \lm}\\
& = & e^{- n \; \vep_a} \\
& < & \exp( n \; g (\vep_a, \vep_a) )\\
& < & \exp \li ( n \; g \li ( \vep_a, \f{\vep_a}{\vep_r} \ri ) \ri ) < \f{\de}{2} \eee where the second inequality follows from Lemma
(\ref{lem3}).   Hence, \[ \Pr \li \{  \li | \wh{\bs{\lm}}  - \lm \ri | \geq \vep_a \; \& \; \li | \wh{\bs{\lm}}  - \lm \ri | \geq \vep_r \lm \ri
\} < \f{\de}{2} < \de.
\]

\bsk

In Case (iii), by Lemma (\ref{lem2}), Lemma (\ref{lem3})  and Lemma (\ref{lem4}), we have \bee \Pr \li \{  \li | \wh{\bs{\lm}}  - \lm \ri | \geq
\vep_a \; \& \; \li | \wh{\bs{\lm}} - \lm \ri | \geq \vep_r \lm
\ri \} & = & \Pr \{ \wh{\bs{\lm}} \leq \lm -  \vep_a \} + \Pr \{ \wh{\bs{\lm}} \geq \lm  + \vep_a \}\\
& \leq & \exp( n \; g (- \vep_a, \lm) ) + \exp( n \; g (\vep_a, \lm) )\\
& < & \exp \li ( n \; g \li (- \vep_a, \f{\vep_a}{\vep_r} \ri ) \ri ) + \exp \li ( n \; g \li ( \vep_a, \f{\vep_a}{\vep_r} \ri ) \ri )\\
& < & 2 \exp \li ( n \; g \li ( \vep_a, \f{\vep_a}{\vep_r} \ri ) \ri ) < \de. \eee

In Case (iv), by Lemma (\ref{lem5}), Lemma (\ref{lem6})  and Lemma (\ref{lem4}), we have  \bee \Pr \li \{  \li | \wh{\bs{\lm}}  - \lm \ri | \geq
\vep_a \; \& \; \li | \wh{\bs{\lm}}  - \lm \ri | \geq \vep_r \lm
\ri \} & = & \Pr \li \{  \li | \wh{\bs{\lm}}  - \lm \ri | \geq \vep_r \lm  \ri \}\\
& = & \Pr \{ \wh{\bs{\lm}} \leq (1 - \vep_r) \lm  \} + \Pr \{ \wh{\bs{\lm}} \geq (1 + \vep_r) \lm  \}\\
& \leq & \exp( n \; g (- \vep_r \lm, \lm) ) + \exp( n \; g (\vep_r \lm, \lm) )\\
& < & \exp \li ( n \; g \li (- \vep_a, \f{\vep_a}{\vep_r} \ri ) \ri ) + \exp \li ( n \; g \li ( \vep_a, \f{\vep_a}{\vep_r} \ri ) \ri )\\
& < & 2 \exp \li ( n \; g \li ( \vep_a, \f{\vep_a}{\vep_r} \ri ) \ri ) < \de. \eee

Therefore, we have shown $\Pr \li \{  \li | \wh{\bs{\lm}}  - \lm \ri | \geq \vep_a \; \& \;  \li | \wh{\bs{\lm}}  - \lm \ri | \geq \vep_r \lm
\ri \} < \de$ for all cases. This completes the proof of Theorem 1.

\end{document}